\documentclass[11pt,article]{article}

\renewcommand{\theequation}{\thesection.\arabic{equation}}

\usepackage{amsmath}
\usepackage{amsthm}
  \usepackage{paralist}
  \usepackage{graphics} 
  \usepackage{epsfig} 
\usepackage{graphicx}
\usepackage{caption}
\usepackage{subcaption}
\usepackage{epstopdf}
 \usepackage[colorlinks=true]{hyperref}
 \usepackage{multirow}
\input{amssym.tex}

\newtheorem{theorem}{Theorem}[section]

\newtheorem{proposition}{Proposition}

 \numberwithin{equation}{section}

\newcommand{\keywords}

\def\bc{\begin{center}}       \def\ec{\end{center}}
\def\ba{\begin{array}}        \def\ea{\end{array}}
\def\be{\begin{equation}}     \def\ee{\end{equation}}
\def\bea{\begin{eqnarray}}    \def\eea{\end{eqnarray}}
\def\beaa{\begin{eqnarray*}}  \def\eeaa{\end{eqnarray*}}

\def\mathbb{\Bbb}

\begin{document}

\title{\bf Pattern formation in Keller--Segel chemotaxis models with logistic growth\footnote{Published in IJBC, Feb, 2016.  Updated}}
\author{LING JIN\footnote{lingjin@2011.swufe.edu.cn. LJ is currently at Department of Mathematics, University of Oklahoma, Norman;}, QI WANG\footnote{qwang@swufe.edu.cn, the corresponding author.  QW is partially supported by NSF-China (Grant 11501460), the Project sponsored by SRF for ROCS, SEM and the Project (No.15ZA0382) from Department of Education, Sichuan China;}, ZENGYAN ZHANG\footnote{zyz@2011.swufe.edu.cn. ZYZ is currently at Department of Applied Mathematics, Hong Kong Polytech University, Hung Hom, Hong Kong.}
\\
Department of Mathematics\\
Southwestern University of Finance and Economics\\
555 Liutai Ave, Wenjiang, Chengdu, Sichuan 611130, China
}

\date{}
\maketitle

\abstract
In this paper we investigate pattern formation in Keller--Segel chemotaxis models over a multi--dimensional bounded domain subject to homogeneous Neumann boundary conditions.  It is shown that the positive homogeneous steady state loses its stability as chemoattraction rate $\chi$ increases.  Then using Crandall--Rabinowitz local theory with $\chi$ being the bifurcation parameter, we obtain the existence of nonhomogeneous steady states of the system which bifurcate from this homogeneous steady state.  Stability of the bifurcating solutions is also established through rigorous and detailed calculations.  Our results provide a selection mechanism of stable wavemode which states that the only stable bifurcation branch must have a wavemode number that minimizes the bifurcation value.  Finally we perform extensive numerical simulations on the formation of stable steady states with striking structures such as boundary spikes, interior spikes, stripes, etc.  These nontrivial patterns can model cellular aggregation that develop through chemotactic movements in biological systems.

\textbf{Keywords: pattern formation, steady state, bifurcation, chemotaxis model, logistic growth}

\section{Introduction}
In numerous biological and ecological contexts, self--organized patterns can emerge and develop through cellular chemotactic movements.  During chemotaxis,  bacteria or cellular organisms sense the chemicals in their environment and orient their migration along the concentration gradient of the chemical.  Such directed movements are observed in many other living microorganisms.  For example, during the aggregation phase of slime molds, the amoebae are attracted by cyclic AMP, a chemical produced by the cells and consumed by certain enzymes in the environment.  Cellular chemotactic movements have also been observed in a wide range of other phenomena, such as wound healing, embryonic development and cancer growth of tumour cells.  In 1970s, Keller and Segel [1971] initiated the mathematical modelling of chemotaxis and studied the traveling bands of E. coli.  Since then, Keller--Segel chemotaxis model has attracted significant interest from numerous scholars over the past few decades and a variety of modifications have been proposed in light of the mathematical relevance to the biology.  We refer to the well--known survey papers [Horstmann, 2003] and [Hillen \& Painter, 2009] as well as the references therein for detailed discussions on chemotaxis models.  In this paper, we study the following Keller--Segel chemotaxis model with logistic cellular growth
\begin{equation}\label{11}
\left\{
\begin{array}{ll}
u_t=\nabla \cdot(D_1 \nabla u-\chi \phi(u,v) \nabla  v)+\mu u(\bar{u}-u),&x \in \Omega,~t>0,     \\
v_t= D_2 \Delta v-\alpha v+ f(u),&x \in \Omega,~t>0,\\
\frac{\partial u}{\partial \textbf{n}}=\frac{\partial v}{\partial \textbf{n}}=0,& x \in \partial \Omega,~t>0,\\
u(x,0)=u_0(x) \geq 0,~v(x,0)=v_0(x) \geq 0,&  x\in \Omega,
\end{array}
\right.
\end{equation}
where $\Delta=\sum_{i=1}^N \partial^2 \slash \partial x_i^2$ is the Laplace operator; $\nabla=(\partial \slash\partial x_1,...,\partial \slash\partial x_N)$ is the gradient operator; $\Omega \subset \mathbb{R}^N$, $N\geq1$, is a bounded domain with piecewise smooth boundary $\partial \Omega$.  $u(x,t)$ and $v(x,t)$ denote the cell population density and chemical concentration at space--time $(x,t)$ respectively.  $D_1>0$ is the so--called cellular motility and it measures the tendency of cells to move
randomly over the habitat.  $D_2>0$ is the diffusion rate of the chemical; the constant $\chi$ measures the attraction rate of the chemical gradient on the directed cellular movement; moreover $\chi>0$ if the chemical is attractive to the cells and $\chi<0$ if the chemical is repulsive.  We will focus on chemoattractive chemical and assume that $\chi>0$ throughout this paper.  $\phi(u,v)$ is called the sensitivity function and it reflects the variation of cellular sensitivity with respect to the cell population density and levels of chemical concentration.  The positive constants $\mu$ and $\bar u$ measure the intrinsic growth rate and the carrying capacity of the habitat for cells respectively; $\alpha>0$ models the phenomenon that the chemical is consumed by certain enzymes in the environment; $f(u)$ is the creation and degradation rate of the chemical with respect to the cell population.  For the sake of mathematical modeling, we assume that the function $\phi$ is $C^3$--smooth and $f$ is $C^2$--smooth such that
\[\phi(u,v)>0,f(u)>0, \text{~for~}u>0,v>0.\]
The non--flux boundary conditions in (\ref{11}) represent an enclosed environment.

One of the most interesting phenomena in chemotaxis is the cellular aggregation in which spatially homogeneously distributed cells move to combine with each other and eventually form a stable concentrated fruiting body.  This can be modeled by the spatial--temporal pattern formation in the time--dependent system (\ref{11}), i.e., solutions $(u(x,t),v(x,t))$ exist globally in time and converge to bounded nonconstant steady states as $t \rightarrow \infty$.  Positive steady states with spiky or concentration structures can be used to model the aggregation of cells.

For $N=1$, it is well known that the solution to (\ref{11}) exists globally and is uniformly bounded in time--see [Osaki \& Yagi, 2001] for example.  For $N=2$, Osaki, \emph{et al}. [2002] proved the global existence of (\ref{11}); moreover, they obtained a globally exponential attractor when the sensitivity function $\phi(u,v)$ is smooth and has uniformly bounded derivatives up to the second order.  For $N\geq3$, Winkler [2010a] established the unique global solutions (\ref{11}) for all smooth initial data if $\mu$ is sufficiently large.  See [Horstmann, 2011] and the references therein for more works.  In loose terms, the literature suggests that blow--ups (finite time or infinite time) in chemotaxis systems can be inhibited by suitable degradation in the cellular kinetics, however, whether or not it is sufficient when chemo--attraction rate $\chi$ is large is not completely understood, in particular over high space dimensions.  For example, Winkler [2011] studied a parabolic--elliptic system similar as (\ref{11}), with $\mu u(\bar u -u)$ being replaced by a term $\mu u(\bar u-u^\kappa)$, $-\alpha v$ by $m(t)=\frac{\int_\Omega u(x,t)dx}{\vert \Omega\vert}$ and $f(u)=u$.  It is shown that when $N\geq 5$, $\kappa<\frac{3}{2}+\frac{1}{2N-2}$ and $v_t=0$, (\ref{11}) has smooth local--in--time solution that blows up in finite time for properly imposed initial data.  Global existence or blow--ups of (\ref{11}) for all parameters remains open so far and we refer to [Xiang, 2015] for some new results.

In contrast to the global existence of (\ref{11}), it is the main purpose of our paper to study the existence and stability of its nonconstant positive stationary solutions which satisfy the following equations
\begin{equation}\label{12}
\left\{
\begin{array}{ll}
\nabla \cdot(D_1 \nabla u-\chi \phi(u,v) \nabla  v)+\mu u(\bar{u}-u)=0,&x \in \Omega,    \\
D_2 \Delta v-\alpha v+ f(u)=0,&x \in \Omega, \\
\frac{\partial u}{\partial \textbf{n}}=\frac{\partial v}{\partial \textbf{n}}=0,& x \in \partial \Omega.
\end{array}
\right.
\end{equation}
In particular, we are interested in the effect of chemotaxis rate on the existence of its nonconstant positive solutions.  System (\ref{12}) has been studied by various authors.  For $D_1=D_2=\alpha=1$, $\phi(u,v)=f(u)=u$ and, Tello and Winkler [2007] obtained infinitely many branches of local bifurcation solutions to (\ref{12}) with $\mu>0$ if $N\leq 4$ and with $\mu>\frac{N-4}{N-2}\chi$ if $N>4$.  For $\phi(u,v)=f(u)=u$, Kuto \emph{et al}. [2012] constructed local bifurcation branches of strip and hexagonal steady states when the domain $\Omega$ is a rectangle in $\mathbb R^2$.  Ma \emph{et al}. [2012] studied the model with a volume-filling effect with $\phi(u,v)=u(1-u)$ and $f(u)=\beta u$, $\beta$ being a positive constant.  They applied degree method to obtain nonconstant positive steady states and established a selection mechanism of the wave modes for $\Omega$ being an interval in $\mathbb{R}^1$.  They also showed that the nontrivial solution is stable only at the first branch and its principal wave mode must be a positive integer that minimizes the bifurcation parameter $\chi$.  The existence and stability of nonconstant positive solutions of (\ref{12}) over one--dimensional finite domain is rigorously investigated in [Wang \emph{et al}., 2015] through bifurcation analysis.  In [Tsujikawa \emph{et al.} 2015], the authors studied a shadow system of (\ref{12}) in the limiting case that a diffusion coefficient and chemotactic intensity grow to infinity through bifurcation, singular perturbation and a level set analysis.  In this paper, we study the existence and stability of nonconstant positive solutions of (\ref{12}) over general domain $\Omega\subset \mathbb R^N$, $N\geq2$.  Compared to the previous works, our work is featured with the rigorous stability analysis  of the nontrivial patterns over general multi--dimensional domains.

Our paper is organized as follows.  In section \ref{section2}, we establish the existence of nonconstant positive solutions of (\ref{12}) by local bifurcation theory--see Theorem \ref{theorem21}.  It is shown that the positive homogeneous solution loses its stability as the chemo--attraction rate passes a threshold value.  Our results indicate that chemotaxis drives the formation of positive nontrivial patterns of the system.  In section \ref{section3}, we proceed to investigate stability of the bifurcating solutions--see Theorem \ref{theorem31}.  Section \ref{section4} is devoted to numerical simulations of pattern formations in the time--dependent system (\ref{11}).  Stable steady states with striking structures such as spikes, stripes etc. are obtained to model the aggregations of cells.  Finally, we include some concluding remarks and propose questions for future studies in section \ref{section5}.

\section{Existence of nonconstant positive steady states}\label{section2}

Pattern formation is a ubiquitous phenomenon observed in natural sciences such as chemistry, physics and biology, etc.  Diffusion--driven instability in reaction--diffusion systems is widely accepted as a mechanism for spatial patterns throughout a wide range of biological and ecological systems--see [Benson \emph{et al.}, 1989], [Okubo, 1980], [Segel \& Jackson, 1977] and the references therein.  All these works are based on Turing's 1952 seminal work "the chemical basis of morphogenesis" [Turing, 1952].  Turing's pioneering idea was that diffusion, which is a smoothing process for single equation, can destabilize constant equilibrium through chemical reactions, then spatially inhomogeneous solutions may emerge and form stable nontrivial patterns.  The phenomenon of Turing's instability (diffusion--driven instability) has been observed in the Gierer--Meinhardt model, Lengyel-Epstein model, as well as diffusive predator--prey models.

System (\ref{12}) has a unique positive constant solution
\[(\bar u,\bar v)=(\bar u,f(\bar u)/\alpha)\]
and we shall see that the instability of this equilibrium is driven by the chemotactic attraction.  To this end, we first show that no Turing's instability occurs for (\ref{12}) with $\chi=0$.  Actually, it is shown in [Ma \emph{et al}., 2013] that $(\bar u,\bar v)$ is a globally attractor of (\ref{11}) if $\chi=0$, therefore, system (\ref{12}) without chemotaxis does not have any nonconstant steady state.  Moreover, this conclusion also holds for $\chi>0$ being small from standard dynamics theory.  For example, similar to the arguments in Theorem 4.3 of [Tello \& Winkler, 2007], one can show that there exists a positive number $\chi_*$ such that (\ref{11}) has no nonconstant positive solution if $\chi \in (0,\chi_*)$.

\subsection{Linearized stability analysis of the homogeneous solution $(\bar u,\bar v)$}

It is well known that random movements (diffusions) tend to stabilize spatially homogeneous solutions for single equations and for reaction--diffusion systems when there is no Turing's instability, however directed movements (chemotaxis) have the effect of destabilizing the homogeneous solutions in general.  Then spatially inhomogeneous solutions may arise through bifurcation as the homogeneous one becomes unstable.

To study the regime of $\chi$ under which spatial patterns of (\ref{11}) arises, we shall perform bifurcation analysis on (\ref{12}) by taking $\chi$ as the parameter and to this end we start with the instability analysis of the homogeneous state $(\bar u,\bar v)$.  First of all, we recall the following Neumann eigenvalue problem
\begin{equation}\label{21}
\left\{
\begin{array}{ll}
-\Delta \Phi =\lambda \Phi,&x \in \Omega,     \\
\frac{\partial \Phi}{\partial \textbf{n}}=0,& x \in \partial \Omega.
\end{array}
\right.
\end{equation}
It is well known that the Neumann Laplacian has a discrete spectrum of infinitely many non-negative eigenvalues which form a strictly increasing sequence
\[0=\lambda_0<\lambda_1<\lambda_2<...<\lambda_k<...\rightarrow \infty;\]
moreover, the eigenvalues are given by the variational structures
\[\lambda_0=\inf_{u\in H^1,\int_\Omega u\neq 0}\frac{\int_\Omega \vert \nabla u\vert ^2dx}{\int_\Omega u^2 dx}\]
and
\[\lambda_k=\inf_{u\in H^1,u\perp \Phi_i,i=1,..k-1}\frac{\int_\Omega \vert \nabla u\vert ^2dx}{\int_\Omega u^2 dx},\]
where $\Phi_i$ is an eigenfunction corresponding to $\lambda_i$; furthermore, $\{\Phi_k\}_{k=1}^\infty$ form a complete orthogonal basis of $L^2(\Omega)$.  Through the rest of this paper, we denote $\{(\lambda_k,\Phi_k)\}_{k=0}^\infty$ as eigen--pairs to (\ref{21}) such that $\lambda_k$ is simple and $\Phi_k$ is normalized with
\[\int_\Omega \Phi_k^2 dx=1,k\in \mathbb N^+.\]
Let $(u,v)=(\bar u,\bar v)+(U,V)$, where $U$ and $V$ are spatially inhomogeneous perturbations of $(\bar u,\bar v)$ in the $H^2(\Omega)$ norm, then $(U,V)$ satisfies
\begin{equation*}
\left\{
\begin{array}{ll}
U_t \approx D_1\Delta U-\chi\phi(\bar u,\bar v)\Delta V-\mu \bar u U,&x \in \Omega, t>0,     \\
V_t \approx D_2\Delta V-\alpha V+f'(\bar u)U,& x \in \Omega, t>0,\\
\frac{\partial U}{\partial \textbf{n}}=\frac{\partial V}{\partial \textbf{n}}=0,& x \in \partial \Omega, t>0.
\end{array}
\right.
\end{equation*}
We have the following results on the linearized instability of $(\bar u, \bar v)$.  Here the stability or instability refers that of the homogeneous solution viewed as an equilibrium of (\ref{11}).
\begin{proposition}\label{proposition1}
 Suppose that $\phi(\bar u,\bar v)>0$ and $f'(\bar u)\neq 0$.  Then the homogeneous solution $(\bar u,\bar v)$ is unstable if and only if
\begin{equation}\label{22}
\chi>\chi_0=\min_{k \in \mathbb{N^+} } \frac{(D_1\lambda_k+\mu \bar u)(D_2\lambda_k+\alpha)}{f'(\bar u)\phi(\bar u,\bar v)\lambda_k}.
\end{equation}
\end{proposition}
\begin{proof}
According to the standard linearized stability analysis--see [Simonett, 1975] e.g., the stability of $(\bar u,\bar v)$ is determined by the eigenvalues of the following matrix,
\begin{equation}\label{23}
 \mathcal H_k=
\begin{pmatrix} -D_1\lambda_k-\mu\bar u& \chi\phi(\bar u, \bar v)\lambda_k\\
f'(\bar u)&-D_2\lambda_k-\alpha \end{pmatrix},
\end{equation}
where $\lambda_k$ is the $k$--th eigenvalue of $-\Delta$ over $\Omega$ under the Neumann boundary conditions.  In particular, $(\bar u,\bar v)$ is stable if eigenvalue of each $\mathcal H_k$ has negative real part and is unstable if $\mathcal H_{k}$ has an eigenvalue with positive real part for some $k \in \mathbb{N}^+$.  The characteristic polynomial of (\ref{23}) takes the form
\[p_k(\lambda)=\lambda^2-\text{Tr}_k\lambda+\text{Det}_k,\]
where \[\text{Tr}_k=(D_1+D_2)\lambda_k+\mu \bar u+\alpha >0,\]
and \[\text{Det}_k=(D_1\lambda_k+\mu \bar u)(D_2\lambda_k+\alpha)-\chi\phi(\bar u,\bar v)f'(\bar u)\lambda_k,\]
therefore $p_k(\lambda)$ has a positive root if and only if $\text{Det}_k<0$.  Then (\ref{22}) readily follows through simple calculations and this finishes the proof of the proposition.
\end{proof}
Proposition \ref{proposition1} states that $(\bar u,\bar v)$ changes its stability as $\chi$ crosses $\chi_0$.  It also indicates that chemo--attraction can destabilizes the homogeneous steady state while chemorepulsion (i.e. when $\chi<0$) stabilizes the constant equilibrium.  We also note that $\chi_0$ given by (\ref{22}) is always positive if $f'(\bar u)>0$.

\subsection{Bifurcation of steady states}
As we have shown above, chemotactic movement has the effect of destabilizing $(\bar u,\bar v)$ which becomes unstable as $\chi$ surpasses the minimum value $\chi_0$.  The linearized instability of $(\bar u,\bar v)$ in (\ref{11}) is insufficient to guarantee the existence of spatially inhomogeneous steady states.  We are concerned with the formation of stable spatially inhomogeneous steady states of (\ref{11}) through bifurcations as $\chi$ increases.  Clearly, the emergence of such spatially inhomogeneous solution is due to the effect of large chemotaxis rate $\chi$.  We refer this as chemotaxis--induced pattern in the sense of Turing's instability.

In this section, we apply the local bifurcation theory in [Crandall \& Rabinowitz, 1971] to seek non-constant positive solutions to the stationary reaction--advection--diffusion system of (\ref{12}).  In order to carry out the bifurcation analysis, we first introduce Sobolev space
\begin{equation}\label{24}
\mathcal{X}=\{w \in W^{2,p}(\Omega) ~\vert  \partial_\textbf{n}w=0, x \in \partial \Omega\},
\end{equation}
where $p>N$.  Taking $\chi$ as the bifurcation parameter, we rewrite (\ref{12}) into the following abstract form
\[\mathcal{F}(u,v,\chi)=0,~(u,v,\chi) \in \mathcal{X}  \times \mathcal{X} \times \mathbb{R},\]
where
\begin{equation}\label{25}
\mathcal{F}(u,v,\chi) =\left(
 \begin{array}{c}
\nabla\cdot(D_1\nabla u-\chi \phi(u,v) \nabla v)+\mu u(\bar u-u)\\
D_2\Delta v-\alpha v+f(u)
 \end{array}
 \right).
 \end{equation}
It is easy to see that $\mathcal{F}(\bar u,\bar v,\chi)=0$ for any $\chi \in \mathbb{R}$ and $\mathcal{F}:\mathcal{X}\times \mathcal{X} \times \mathbb{R} \rightarrow L^p(\Omega)\times L^p(\Omega)$ is analytic.  For any fixed $(\hat u,\hat{v})\in\mathcal{X}\times\mathcal{X}$, we have from straightforward calculations that the Fr\'echet derivative of $\mathcal{F}$ is given by
\begin{equation}\label{26}
D_{(u,v)}\mathcal{F}(\hat{u},\hat{v},\chi)(u,v)=\left(
\begin{array}{c}
\nabla \cdot\Big(D_1 \nabla u-\chi\big((\phi_u(\hat{u},\hat{v})u+\phi_v(\hat{u},\hat{v})v\big)\nabla \hat{v}\\+\phi(\hat{u},\hat{v})\nabla v\big)\Big)+\mu(\bar{u}-2\hat{u})u\\
D_2\Delta v-\alpha v+f'(\hat u)u
\end{array}
\right);
\end{equation}
moreover it is easy to see that $D_{(u,v)}\mathcal F(\hat u,\hat v, \lambda)$ is continuously differentiable with respect to $(u,v)$ and $\chi$ in $\mathcal{X}\times \mathcal{X} \times \mathbb R$.

To apply the Crandall--Rabinowitz bifurcation theory, we need to verify that $D_{(u,v)}\mathcal{F}(\bar u,\bar v,\chi)(u,v)$ is Fredholm with index 0 and it has zero as a simple eigenvalue.  Before showing this, we want to remark that the arguments developed by [Shi \& Wang, 2009] make the original bifurcation theory applicable in the Sobolev spaces above (and also in H\"older spaces thanks to elliptic Sobolev embeddings).

For local bifurcations to occur at $(\bar u,\bar v,\chi)$, the following nontriviality condition is needed for the null space of $D_{(u,v)}\mathcal{F}(\bar{u},\bar{v},\chi)$ to hold
\[\mathcal{N}\big(D_{(u,v)}\mathcal{F}(\bar{u},\bar{v},\chi)\big) \neq \{(0,0)\}.\]
To verify this condition, we first take $(\hat u, \hat v)=(\bar u,\bar v)$ in (\ref{26}) and have that Fr\'echet derivative of $\mathcal{F}$ there reads
\begin{equation}\label{27}
D_{(u,v)}\mathcal{F}(\bar u,\bar v,\chi)(u,v)=\left(
\begin{array}{c}
D_1\Delta u-\chi\phi(\bar u, \bar v)\Delta v-\mu \bar u u \\
D_2\Delta v-\alpha v+f'(\bar u)u
\end{array}
\right),
\end{equation}
hence the nontriviality condition implies that there exists some nontrivial solutions $(u,v)$ to the following system
\begin{equation}\label{28}
\left\{
\begin{array}{ll}
D_1\Delta u -\chi\phi(\bar{u},\bar{v})\Delta v-\mu\bar{u}u=0,& x\in\Omega,\\
D_2  \Delta v-\alpha v+f'(\bar u)u=0,& x\in\Omega,\\
\frac{\partial u}{\partial \textbf{n}}=\frac{\partial v}{\partial \textbf{n}}=0,& x \in \partial \Omega.
\end{array}
\right.
\end{equation}
We expand the $u$ and $v$ into their eigen--expansions
\[u(x)=\sum_{k=1}^{\infty} t_{k}\Phi_k,v(x)=\sum_{k=1}^{\infty} s_{k}\Phi_k,\]
where $t_k$ and $s_k$ are constants to be determined.  Substitute the series into (\ref{28}) and we obtain
\begin{equation}\label{29}
\begin{pmatrix}
-D_1 \lambda_k-\mu \bar u & \lambda_k\chi\phi(\bar u, \bar v)   \\
~~\\
f'(\bar u) & -D_2\lambda_k-\alpha
\end{pmatrix}
\begin{pmatrix}
t_{k}\\
~~\\
s_{k}
\end{pmatrix}=\begin{pmatrix}
0\\
~~\\
0
\end{pmatrix},
\end{equation}
where $\lambda_k$ is the $k$--th eigenvalue of (\ref{21}).  $k=0$ can be easily ruled out if (\ref{27}) admits nonzero solutions since $\lambda_0=0$.  For each $k \in \mathbb{N}^+$, the coefficient matrix in (\ref{29}) is singular and one has that
\begin{equation}\label{210}
\chi=\bar \chi_{k}=\frac{(D_1\lambda_k+\mu \bar u)(D_2\lambda_k+\alpha)}{f'(\bar u)\phi(\bar u,\bar v)\lambda_k}, k\in \mathbb N^+;
\end{equation}
moreover the null space $\mathcal{N}(D_{(u,v)}\mathcal{F}(\bar{u},\bar{v},\bar \chi_{k}))$ is one--dimensional and it has a span \[\mathcal{N}(D_{(u,v)}\mathcal{F}(\bar{u},\bar{v},\bar \chi_{k}))=\text{span}\left\{ (\bar{u}_{k}, \bar{v}_{k}) \right\},\]
where
\begin{equation}\label{211}
\bar{u}_{k}=Q_k\Phi_k,~\bar{v}_{k}=\Phi_k,
\end{equation}
with $\Phi_k$ being the $k$--th Neumann Laplace eigenfunction and
\begin{equation}\label{212}
Q_{k}=\frac{D_2\lambda_k+\alpha}{f'(\bar u)}, k \in \mathbb{N}^+.
\end{equation}
Now it follows from Corollary 2.11 in [Shi \& Wang, 2009] that $D_{(u,v)}\mathcal{F}(\hat u,\hat{v},\chi)$ is a Fredholm operator with index 0--see also the statements of Theorem 3.3 and Remark 3.4 in that paper.  Therefore we have that codim $\mathcal R(D_{(u,v)}\mathcal{F}(\hat u,\hat{v},\chi))=1$.  Having the potential bifurcation values in (\ref{210}), we now verify in the following theorem that local bifurcation does occur at $(\bar u,\bar v,\bar \chi_{k})$ for each $k \in \mathbb{N^+}$, which establishes nonconstant positive solutions to (\ref{12}).
\begin{theorem}\label{theorem21}
Let $\lambda_i$ be the $i$--th simple eigenvalue of (\ref{21}) and assume that $D_1D_2\lambda_i\lambda_j\neq \alpha \mu \bar u$ for all $i\neq j\in \mathbb N^+$.  Then the solutions of (\ref{12}) around $(\bar u,\bar v,\bar \chi_k)$ consists precisely of the continuous curve $\Gamma_{k}(s)=(u_{k}(s), v_{k}(s), \chi_{k}(s))$, $s\in(-\delta,\delta)$, where
\begin{equation}\label{213}
\chi_{k}(s)=\bar \chi_{k}+o(s),~(u_{k}(s,x),v_{k}(s,x))=(\bar u,\bar v)+s(Q_{k},1)\Phi_k+o(s),
\end{equation}
and $(u_{k}(s,x),v_{k}(s,x))-(\bar u,\bar v)-s(Q_{k},1)\Phi_k$ in the closed complement of the null space $\mathcal{N}(D_{(u,v)}\mathcal{F}(\bar{u},\bar{v},\bar \chi_{k}))$ which is explicitly given by
\begin{equation}\label{214}
\mathcal{Z}=\big\{(u,v)\in \mathcal{X} \times \mathcal{X}~ \big \vert \int_{\Omega} u \bar u_{k}+v\bar v_{k} dx=0\big\}
\end{equation}
with $(\bar u_{k},\bar v_{k})$ given by (\ref{211}).
\end{theorem}
\begin{proof}
In order to apply the local theory of Crandall and Rabinowitz [1971], we have verified all but the following transversality condition,
\begin{equation}\label{215}
\frac{d}{d \chi} \left(D_{(u,v)}\mathcal{F}(\bar{u},\bar{v},\chi)\right)(\bar u_{k},\bar v_{k})\vert_{\chi=\bar \chi_{k}} \notin \mathcal{R}(D_{(u,v)}\mathcal{F}(\bar{u},\bar{v},\bar \chi_{k})).
\end{equation}
We argue by contradiction.  If not, there must exist a nontrivial solution $(\tilde{u},\tilde{v})$ to the following problem
\begin{equation}\label{216}
\left\{
\begin{array}{ll}
D_1 \Delta \tilde{u}-\bar \chi_{k}\phi(\bar{u},\bar{v})\Delta \tilde{v}-\mu \bar{u}\tilde{u}
=- \lambda_k\phi(\bar{u},\bar{v}) \Phi_k ,&x \in\Omega,\\
D_2  \Delta\tilde{v}- \alpha \tilde v+f'(\bar u)\tilde u=0,&x \in\Omega,\\
\frac{\partial \tilde u}{\partial \textbf{n}}=\frac{\partial \tilde{v}}{\partial \textbf{n}}=0,&x \in \partial \Omega.
\end{array}
\right.
\end{equation}
Multiplying both equations in (\ref{216}) by $\Phi_k$ and then integrating them over $\Omega$ by parts, we have that
\begin{equation}\label{217}
\begin{pmatrix}
-D_1 \lambda_k-\mu\bar{u} & \lambda_k\bar \chi_{k}\phi(\bar{u},\bar{v})  \\
~~\\
f'(\bar u) &-D_2\lambda_k-\alpha
\end{pmatrix}
\begin{pmatrix}
\int_\Omega \tilde u \Phi_kdx \\
~~\\
\int_\Omega \tilde v \Phi_kdx
\end{pmatrix}=\begin{pmatrix}
-\lambda_k\phi(\bar u,\bar v) \\
~~\\
0
\end{pmatrix},
\end{equation}where we have used the fact that $\int_\Omega \Phi^2_k dx=1$.  Since $\lambda_k>0$ for $k\in \mathbb N^+$, (\ref{217}) leads to a contradiction in light of (\ref{210}).  This proves (\ref{215}) and the statements of this theorem follows from Theorem 1.7 in [Crandall \& Rabinowitz, 1971].
\end{proof}

We assume that $D_1D_2\lambda_i\lambda_j\neq \alpha \mu \bar u$ for all $i\neq j\in \mathbb N^+$ in order to make sure that $\chi_i\neq \chi_j$.  The bifurcation branch $\Gamma_k(s)$ is local as presented in Theorem \ref{theorem21}.  It is actually a global one according to Theorem 6.1 in [Pejsachowicz \& Rabier, 1998] (See also Remark 4.2 in [Shi \& Wang, 2009]).  From the global bifurcation theory of Crandall and Rabinowitz [1971] and its extended version in [Shi \& Wang, 2009], one can assert that for $\mathcal V$ being an open connected subset of $\mathcal X\times \mathcal X \times \mathbb R$, each curve $\Gamma_k(s)$ is contained in $\mathcal C$ which is a connected component of $\bar {\mathcal S}$, the closure of solution set
$\mathcal S=\{(u,v,\chi)\in\mathcal V: \mathcal F(u,v,\chi)=0,(u,v)\neq(\bar u,\bar v)\}$; moreover, $\mathcal C$ satisfies one of the follows: (i) it is not compact in $\mathcal V$; (ii) it intersects with $(\bar u,\bar v,\chi_*)$ with $\chi_*\neq \bar \chi_k$; or (iii) it contains a point $(\bar u+u,\bar v+v,\chi)$, with $(0,0)\neq (u,v)\in \mathcal Z$.  Here $\mathcal C$ being not compact means that it either intersects $\partial \mathcal V$ or it is unbounded.

As pointed out in [Shi \& Wang, 2009], the global bifurcation results have very useful applications in the study of positive (and monotone) solutions to elliptic PDEs.  However, we want to add that very few mathematical techniques are available in rigorously characterizing these alternatives.  The methodology recently developed by Wang \emph{et al}. in [Chertock \emph{et al}. 2012] and [Xu \& Wang, 2012] turns out to be very useful and user friendly in dealing with certain one--dimensional chemotaxis systems with special structures.  In loose terms, the ideas there can be summarized as follows.  Consider the first bifurcation branch and show that the solutions on the continuum $\mathcal C$ must be monotone by maximum principles.  Therefore the continuum of the first branch must be unbounded.  Then, by deriving a--prior estimates of positive solutions, one can show the continuum $\mathcal C$ extends to infinity and its project onto the $\chi$--axis must be an interval that contains $[\chi_1,\infty)$.  Hence, the first bifurcation branch extends to infinity in the positive direction of $\chi$--axis.  One can then study the behaviors of positive solutions as $\chi \rightarrow \infty$ that has spikes or transition layers.  However, the methodology there does not apply to multi--dimensional domains since the monotonicity argument becomes technically difficult in this case and this is beyond the scope of our paper.  We refer to the recent paper [Ma \& Wang, 2015] on the global bifurcation analysis of (\ref{11}) in 1D.

\section{Stability of bifurcating solutions}\label{section3}

Now we proceed to investigate the stability of the bifurcating solution of (\ref{12}) established in Theorem \ref{theorem21}. Here the stability is that of $(u_k(s,x),v_k(s,x))$ viewed as an equilibrium of system (\ref{11}).  To this end, we first determine the type of each bifurcation branch $\Gamma_{k}(s)$, $k \in \mathbb{N}^+$.  The operator $\mathcal{F}$ defined in (\ref{25}) is $C^4$--smooth if $\phi$ is $C^5$--smooth, hence by Theorem 1.7 from [Crandall \& Rabinowitz 1971], $(u_{k}(s,x),v_{k}(s,x),\chi_{k}(s))$ is a $C^3$--smooth function of $s$ and we have the following asymptotic expansions for $s\in(-\delta,\delta)$:
\begin{equation}\label{31}
\left\{
\begin{array}{ll}
u_{k}(s,x)=\bar{u}+sQ_{k}\Phi_k+s^2\psi_1(x)+s^3\tilde \psi_1(x)+O(s^4),\\
v_{k}(s,x)=\bar{v}+s\Phi_k+s^2\psi_2(x)+s^3\tilde \psi_2(x)+O(s^4),\\
\chi_{k}(s)=\bar \chi_{k}+s\chi'_k(0)+\frac{1}{2}s^2\chi''_k(0)+O(s^2),
\end{array}
\right.
\end{equation}
where $Q_k$ is given by (\ref{212}) and $\Phi_k(x)$ is the $k$--th eigenfunction of the Neumann problem (\ref{21}).  $\psi_i$ and $\tilde \psi_i$, $i=1,2$, belong to $\mathcal{Z}$ as defined in (\ref{214}) and $O(s^4)$ in the first two equations of (\ref{31}) are taken in $W^{2,p}$-topology.  As we shall see later on, the sign of $\chi_k(s)$ around $\bar \chi_k$ determines the stability of the bifurcating solutions.  To this end, we first evaluate $\chi_k'(0)$ and $\chi_k''(0)$ in the following propositions.
\begin{proposition}\label{proposition2}
Let $(\lambda_k,\Phi_k)$ be an eigen--pair of $-\Delta$.  Then $\chi'_k(0)$ is given by
\begin{equation}\label{311}
\chi'_k(0)\lambda_k\phi(\bar u,\bar v)=\Big(\mu Q_k^2-\frac{\lambda_k\bar \chi_k}{2}\big(\phi_u(\bar u,\bar v)Q_k+\phi_v(\bar u,\bar v)\big)\Big)\int_{\Omega}\Phi_k^3 dx.
\end{equation}
\end{proposition}
\begin{proof}
First of all, we have the following facts from Taylor expansions
\begin{equation}\label{32}
\begin{split}
\phi(u,v)=&\phi(\bar u,\bar v)+s\Phi_k\big(\phi_u(\bar u,\bar v)Q_k+\phi_v(\bar u,\bar v)\big)+s^2\Big(\phi_u(\bar u,\bar v)\psi_1+\phi_v(\bar u,\bar v)\psi_2\\
&+\frac{1}{2}\big(\phi_{uu}(\bar u,\bar v)Q_k^2+\phi_{vv}(\bar u,\bar v)+2\phi_{uv}(\bar u,\bar v)Q_k\big)\Phi_k^2\Big)+o(s^2) \\
\end{split}
\end{equation}
and
\begin{equation}\label{33}
\begin{split}
\mu u(\bar u-u)&=-\mu\big(\bar u+sQ_k\Phi_k+s^2\psi_1+O(s^3)\big)\big(sQ_k\Phi_k+s^2\psi_1+O(s^3)\big)  \\
&=s(-\mu \bar u Q_k\Phi_k)+s^2(-\mu \bar u\psi_1-\mu Q_k^2\Phi_k^2)+O(s^3).
\end{split}
\end{equation}

Substituting (\ref{31})--(\ref{33}) into the $u$--equation of (\ref{12}) and equating the $s^2$-terms, we collect that
\begin{equation}\label{34}
\begin{split}
&D_1\Delta\psi_1-\mu\bar u\psi_1-\mu Q_k^2\Phi_k^2+\chi'_k(0)\phi(\bar u,\bar v)\lambda_k\Phi_k\\
=&\bar \chi_k\Big(\phi(\bar u,\bar v)\Delta \psi_2+\big(\phi_u(\bar u,\bar v)Q_k+\phi_v(\bar u,\bar v)\big)\big(\vert\nabla\Phi_k\vert^2-\lambda_k\Phi_k^2\big)\Big).
\end{split}
\end{equation}
Similarly, substituting the expansions into the $v$--equation and equating the $s^2$--terms there, we obtain that
\begin{equation}\label{35}
D_2\Delta\psi_2-\alpha\psi_2+f'(\bar u)\psi_1=0.
\end{equation}

Multiplying (\ref{34}) and (\ref{35}) by $\Phi_k$ and integrating them over $\Omega$ by parts, we have
\begin{equation}\label{36}
\begin{split}
\chi'_k(0)\phi(\bar u,\bar v)\lambda_k=&(D_1\lambda_k+\mu\bar u)\int_\Omega \psi_1\Phi_k dx+\mu Q_k^2\int_{\Omega}\Phi_k^3dx-\lambda_k \bar \chi_k\phi(\bar u,\bar v)\int_\Omega \psi_2\Phi_k dx\\
&+\bar \chi_k\big(\phi_u(\bar u,\bar v)Q_k+\phi_v(\bar u,\bar v)\big)\Big(\int_\Omega \Phi_k\vert \nabla \Phi_k\vert^2 dx-\lambda_k\int_\Omega \Phi_k^3 dx\Big).
\end{split}
\end{equation}
and
\begin{equation}\label{37}
-(D_2\lambda_k+\alpha)\int_{\Omega}\psi_2\Phi_kdx+f'(\bar u)\int_{\Omega}\psi_1\Phi_kdx=0
\end{equation}
respectively; on the other hand, since $(\psi_1,\psi_2)\in \mathcal{Z}$, we have from (\ref{214}) that
\begin{equation}\label{38}
Q_{k}\int_{\Omega} \psi_1 \Phi_kdx+\int_{\Omega} \psi_2 \Phi_kdx=0
\end{equation}
where $Q_{k}=\frac{D_2\lambda_k+\alpha}{f'(\bar u)}$.  It is easy to see that the coefficient matrix of (\ref{37}) and (\ref{38}) is nonsingular from (\ref{212}), therefore we must have that
\begin{equation}\label{39}
\int_{\Omega} \psi_1\Phi_k dx=\int_{\Omega} \psi_2 \Phi_k dx=0, \text{~for~all~} k\in \mathbb{N}^+;
\end{equation}
moreover, it follows from the Neumann eigenvalue problem (\ref{21}) and the straightforward calculations that
\begin{equation}\label{310}
\int_\Omega \Phi_k \vert \nabla \Phi_k\vert^2 dx=\frac{1}{2}\int_\Omega \nabla \Phi_k \nabla\Phi_k^2 dx=-\frac{1}{2}\int_\Omega \Phi_k^2\Delta\Phi_kdx =\frac{\lambda_k}{2}\int_\Omega\Phi_k^3dx;
\end{equation}
then (\ref{311}) follows from (\ref{36}), (\ref{39}) and (\ref{310}).
\end{proof}
Proposition \ref{proposition2} indicates that $\chi'_k(0)=0$ if $\Vert \Phi_k \Vert_{L^3}=0$.  This is the case when $\Omega$ is a one--dimensional finite interval or multi--dimensional rectangle and then the bifurcation is of pitch--fork type hence we need to evaluate $\chi''_k(0)$ for the stability analysis.  We have the following result from straightforward calculations.

\begin{proposition}
If $\chi'_k(0)=0$, then $\chi''_k(0)$ satisfies
\begin{equation}\label{310a}
\begin{aligned}
\frac{1}{2}\phi(\bar u,\bar v)\lambda_k\chi''_k(0)=&\bar{\chi_k}\Big(\phi_u(\bar u,\bar v)Q_k\int_\Omega \psi_2\vert\nabla\Phi_k\vert^2dx-\phi_u(\bar u, \bar v)\int_\Omega \psi_1\vert\nabla\Phi_k\vert^2 dx\\
&-\frac{\lambda_k}{6}\big(\phi_{uu}(\bar u,\bar v)Q_k^2+\phi_{vv}(\bar u, \bar v)+2\phi_{uv}(\bar u, \bar v)Q_k\big)\int_\Omega \Phi_k^4 dx\\
&-\lambda_k\big(\phi_u(\bar u,\bar v)Q_k+\phi_v(\bar u,\bar v)\big)\int_\Omega \psi_2\Phi_k^2 dx\Big)\\
&+2\mu Q_k\int_\Omega \psi_1\Phi_k^2 dx.
\end{aligned}
\end{equation}
\end{proposition}
\begin{proof}
Equating $s^3$ terms in (\ref{12}), together with the fact that $\chi'_k(0)=0$, we arrive at the following system
\begin{equation}\label{310b}
\left\{
\begin{split}
&D_1\Delta\tilde{\psi_1}-\mu\bar u\tilde{\psi_1}-2\mu Q_k\psi_1\Phi_k+\frac{1}{2}\chi''_k(0)\phi(\bar u,\bar v)\lambda_k\Phi_k\\
=&\bar \chi_k\Big(\big(\phi_u(\bar u,\bar v)Q_k+2\phi_v(\bar u,\bar v)\big)\nabla \psi_2\nabla\Phi_k+\phi_u(\bar u,\bar v)\nabla\psi_1\nabla\Phi_k\\
&+\big(\phi_{uu}(\bar u,\bar v)Q_k^2+\phi_{vv}(\bar u, \bar v)+2\phi_{uv}(\bar u, \bar v)Q_k\big)\big(\Phi_k\vert\nabla\Phi_k\vert^2-\frac{1}{2}\lambda_k\Phi_k^3\big)\\
&+\phi(\bar u,\bar v)\Delta \tilde{\psi_2}+\big(\phi_u(\bar u,\bar v)Q_k+\phi_v(\bar u,\bar v)\big)\Delta\psi_2\Phi_k\\
&-\lambda_k\phi_u(\bar u,\bar v)\psi_1\Phi_k-\lambda_k\phi_v(\bar u,\bar v)\psi_2\Phi_k\Big),\\
&D_2\Delta\tilde{\psi_2}-\alpha\tilde{\psi_2}+f'(\bar u)\tilde{\psi_1}=0,\\
&\frac{\partial \tilde{\psi_1}}{\partial \textbf{n}}=\frac{\partial \tilde{\psi_2}}{\partial \textbf{n}}=0.
\end{split}
\right.
\end{equation}
Multiplying the first equation of (\ref{310b}) by $\Phi_k$, we have from the integration by parts that
\begin{equation}\label{310c}
\begin{aligned}
\frac{1}{2}\phi(\bar u,\bar v)\lambda_k\chi''_k(0)=&-D_1 \int_\Omega \Delta\tilde{\psi_1}\Phi_k dx+\mu\bar u\int_\Omega \tilde{\psi_1}\Phi_k dx+2\mu Q_k\int_\Omega \psi_1\Phi_k^2 dx\\
&+\bar \chi_k\Big(\big(\phi_u(\bar u,\bar v)Q_k+2\phi_v(\bar u,\bar v)\big)\int_\Omega \nabla \psi_2\nabla\Phi_k\Phi_k dx+\phi_u(\bar u,\bar v)\int_\Omega \nabla\psi_1\nabla\Phi_k\Phi_k dx\\
&+\big(\phi_{uu}(\bar u,\bar v)Q_k^2+\phi_{vv}(\bar u, \bar v)+2\phi_{uv}(\bar u, \bar v)Q_k\big)\big(\int_\Omega\Phi_k^2\vert\nabla\Phi_k\vert^2 dx-\frac{1}{2}\lambda_k\int_\Omega \Phi_k^4 dx\big)\\
&+\phi(\bar u,\bar v)\int_\Omega \Delta \tilde{\psi_2}\Phi_k dx+\big(\phi_u(\bar u,\bar v)Q_k+\phi_v(\bar u,\bar v)\big)\int_\Omega \Delta\psi_2\Phi_k^2 dx\\
&-\lambda_k\phi_u(\bar u,\bar v)\int_\Omega \psi_1\Phi_k^2 dx-\lambda_k\phi_v(\bar u,\bar v)\int_\Omega \psi_2\Phi_k^2 dx\Big).
\end{aligned}
\end{equation}
By the same arguments that lead to (\ref{39}), we can show that
\[\int_{\Omega} \tilde \psi_1\Phi_k dx=\int_{\Omega} \tilde \psi_2 \Phi_k dx=0, \forall k\in \mathbb{N}^+,\]
and it gives rise to (\ref{310a}) in light of (\ref{310c}).
\end{proof}
In order to evaluate $\chi''_k(0)$, we need to calculate the followings integrals:
\[\int_\Omega \psi_1\Phi_k^2 dx,\int_\Omega \psi_2\vert\nabla\Phi_k\vert^2dx,\int_\Omega \psi_1\vert\nabla\Phi_k\vert^2 dx,\int_\Omega \psi_2\Phi_k^2 dx.\]
Similar as the calculations above, we can have from some straightforward but lengthy calculations
\begin{equation*}
\begin{aligned}
&\begin{pmatrix}
-2D_1 \lambda_k-\mu \bar u & 2D_1 & 2\bar{\chi}_k\lambda_k\phi(\bar u, \bar v) & -2\bar{\chi}_k\phi(\bar u, \bar v)  \\
2D_1\lambda_k^2 & -2D_1\lambda_k-\mu\bar u & -2\bar{\chi}_k\lambda_k^2\phi(\bar u, \bar v) & 2\bar{\chi}_k\lambda_k\phi(\bar u, \bar v)   \\
f'(\bar u) & 0 & -2D_2\lambda_k-\alpha & 2D_2   \\
0 & f'(\bar u) & 2D_2\lambda_k^2 & -2D_2\lambda_k-\alpha
\end{pmatrix}
\begin{pmatrix}
\int_{\Omega}\Phi_k^2\psi_1dx \\
\int_{\Omega}\vert\nabla\Phi_k\vert^2\psi_1dx \\
\int_{\Omega}\Phi_k^2\psi_2dx \\
\int_{\Omega}\vert\nabla\Phi_k\vert^2\psi_2dx
\end{pmatrix}\\
=&\begin{pmatrix}
\Big(\frac{2}{3}\bar{\chi}_k\big(\phi_u(\bar u, \bar v)Q_k+\phi_v(\bar u, \bar v)\big)\lambda_k+\mu Q_k^2\Big) \\
\Big(\frac{2}{3}\lambda_k^2\bar\chi_k\big(\phi_u(\bar u, \bar v)Q_k+\phi_v(\bar u, \bar v)\big)+\frac{1}{3}\lambda_k\mu Q_k^2\Big)\\
0  \\
0
\end{pmatrix}\int_{\Omega}\Phi_k^4dx.
\end{aligned}
\end{equation*}
The calculations for these integrals are straightforward but very complicated, therefore we skip the details here.

Now we are ready to present the following result on the stability of the bifurcating solutions.
\begin{theorem}\label{theorem31}
Suppose that all conditions in Theorem \ref{theorem21} hold.  Assume that $\chi_{k_0}=\min_{k\in\mathbb N^+}\chi_k$, then for all $k\neq k_0$ the bifurcating solution $(u_k(s,x),v_k(s,x))$ is always unstable for $s\in(-\delta,\delta)$.  The bifurcating solution $(u_{k_0}(s,x),v_{k_0}(s,x))$ is asymptotically stable for $s\in(0,\delta)$ and unstable for $s\in(-\delta,0)$ if $\chi_{k_0}'(0)>0$ in (\ref{311}), $(u_{k_0}(s,x),v_{k_0}(s,x))$ is unstable for $s\in(0,\delta)$ and asymptotically stable for $s\in(-\delta,0)$ if $\chi_{k_0}'(0)<0$.  If $\chi'_{k_0}(0)=0$, $(u_{k_0}(s,x),v_{k_0}(s,x))$ is asymptotically stable for $s\in(0,\delta)$ if $\chi''_{k_0}(0)>0$ and unstable if $\chi''_{k_0}(0)<0$.
\end{theorem}
It is easy to see that $\int_{\Omega}\Phi_k^3 dx$ plays an essential role in determining the stability of the bifurcation branch $\Gamma_k(s)$.  Transcritical and pitch bifurcation branches are schematically presented in Figure \ref{fig1} to illustrate the results in Theorem \ref{theorem31}.
\begin{figure}
\minipage{1.\textwidth}\centering
\centering
    \begin{subfigure}[b]{0.3\textwidth}
        \includegraphics[width=\textwidth]{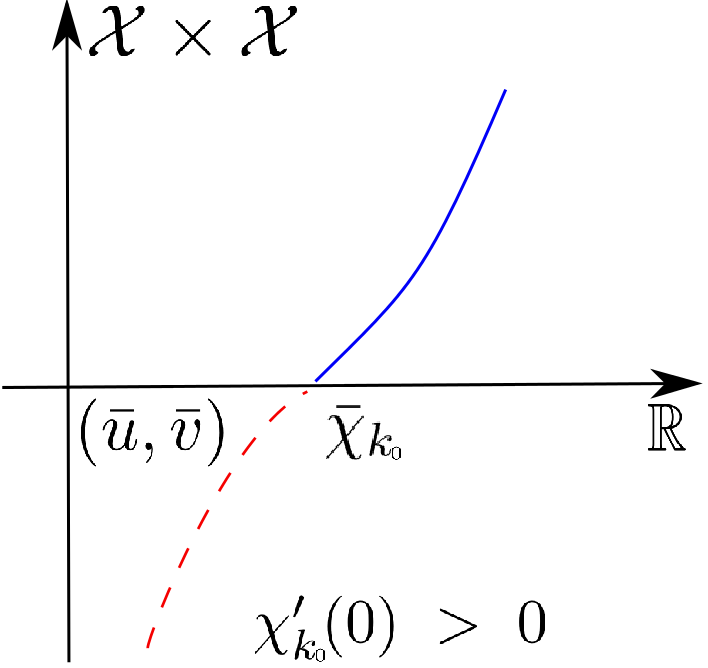}
        \caption{Supercritical bifurcation}
    \end{subfigure}
    \begin{subfigure}[b]{0.3\textwidth}
        \includegraphics[width=\textwidth]{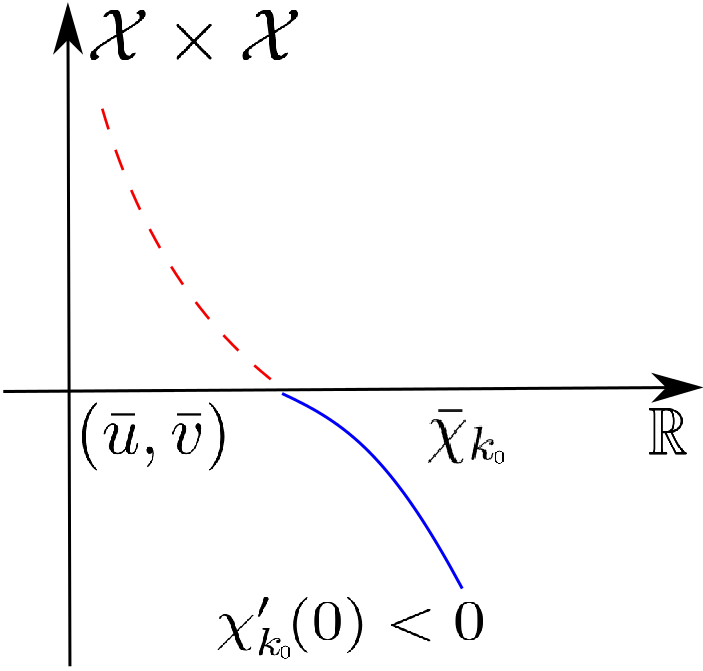}
        \caption{Subcritical bifurcation}
    \end{subfigure}
\endminipage\\
\minipage{1.\textwidth}\centering
\centering
    \begin{subfigure}[b]{0.3\textwidth}
        \includegraphics[width=\textwidth]{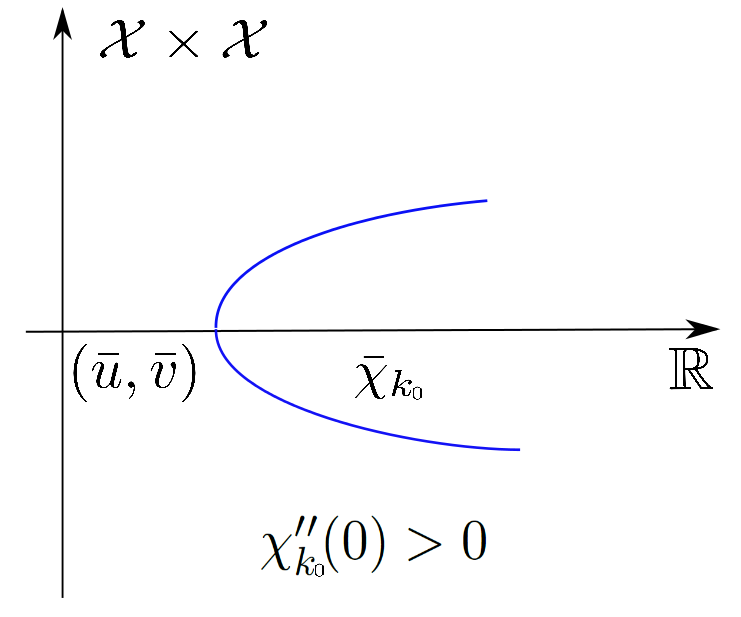}
        \caption{Supercritical bifurcation}
    \end{subfigure}
    \begin{subfigure}[b]{0.3\textwidth}
        \includegraphics[width=\textwidth]{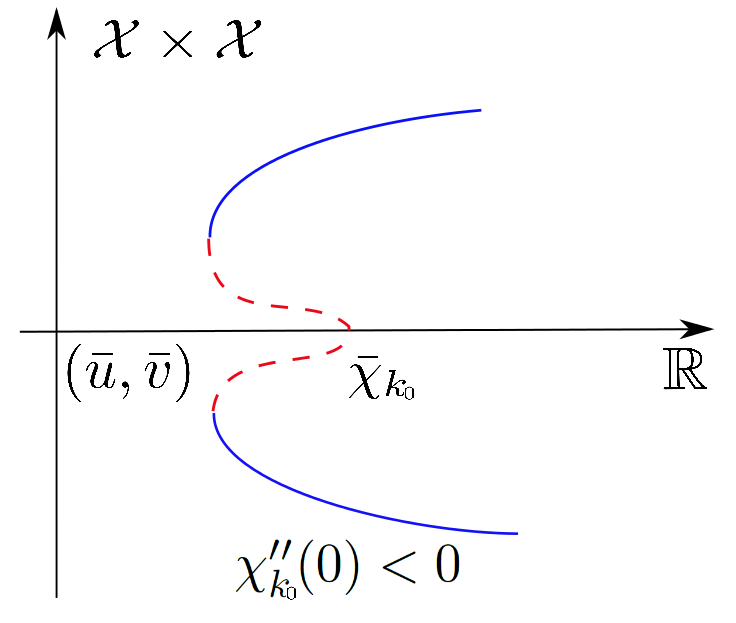}
        \caption{Subcritical bifurcation}
    \end{subfigure}
\caption{Bifurcation branches $\Gamma_{k_0}(s)$ around $(\bar u,\bar v,\bar \chi_{k_0})$.  Transcritical bifurcation branches in (a)--(b) and pitch--fork type bifurcation branches $\Gamma_{k_0}(s)$ in (c)--(d).  The solid curve represents stable bifurcating solution $(u_{k_0}(s,x),v_{k_0}(s,x))$ and the dashed curve represents unstable solution.  Solution on the bifurcation branch changes its stability each time the curve turns its direction.}
\endminipage
\end{figure}\label{fig1}
\begin{proof}Our proof follows that of Theorem 3.2 in [Wang \emph{et al}., 2015].  To prove the instability of solutions $(u_k(s,x),v_k(s,x))$ on bifurcation branch $\Gamma_k(s)$ for $k\neq k_0$, it is equivalent to show that the real part of each eigenvalue of $D_{(u,v)}\mathcal F (u_k(s,x),v_k(s,x),\chi_k(s))$ to the following eigenvalue problem is positive
\begin{equation}\label{312}
D_{(u,v)}\mathcal F (u_k(s,x),v_k(s,x),\chi_k(s))(u,v)=\sigma (u,v),~(u,v)\in \mathcal X \times \mathcal X.
\end{equation}
Indeed, for $s=0$ the linearized stability matrix of (\ref{312}) becomes (\ref{23}) which always has a positive eigenvalue if $k\neq k_0$.  Therefore (\ref{312}) always has a positive eigenvalue for $s$ being small in light of the standard eigenvalue perturbation theory in [Kato, 1996].

Applying the same arguments that lead to the Fredholmness in Theorem \ref{theorem21}, we can show that 0 is a K--simple eigenvalue of $D_{(u,v)}\mathcal F (\bar u,\bar v,\bar \chi_{k_0})$--see Definition 1.2 in [Crandall \& Rabinowitz, 1973].  From Corollary 1.13 in [Crandall \& Rabinowitz, 1973], there exist a small interval $I$ containing $\bar \chi_{k_0}$ and continuously differentiable functions
$(\gamma(\chi),\sigma(s)):I\times(-\delta,\delta)\rightarrow \mathbb R^2$ such that $\gamma=\gamma(\chi)$ is a real eigenvalue of
 \begin{equation}\label{313}
 D_{(u,v)}\mathcal F (\bar u,\bar v,\chi)(u,v)=\gamma(\chi)(u,v),~(u,v)\in \mathcal X \times \mathcal X
 \end{equation}
such that $\gamma(\bar \chi_{k_0})=0$ and $\sigma=\sigma(s)$ is an eigenvalue of (\ref{312}) with $\sigma(0)=0$ ; moreover $\gamma(\chi)$ is the only eigenvalue of (\ref{313}) for any fixed neighbourhood of the origin in the complex plane;  furthermore, the eigenfunction of (\ref{313}) corresponding to $\gamma(\chi)$ depends on $\chi$ smoothly which can be written as $(u(\chi,x),v(\chi,x))$ and is uniquely determined by $(u(\bar \chi_{k_0},x),v(\bar \chi_{k_0},x))=(\bar u_{k_0},\bar v_{k_0})$ and $(u(\chi,x),v(\chi,x))-(\bar u_{k_0},\bar v_{k_0}) \in \mathcal Z$.

According to 1.17 in Theorem 1.16 of [Crandall \& Rabinowitz, 1973], functions $\sigma(s)$ and $-s\chi'_{k_0}(s)\dot{\gamma}(\bar \chi_{k_0})$ have the same zeros and the same sign near $s=0$ and for $\sigma(s) \neq0$,
\[\lim_{s\rightarrow 0}\frac{-s\chi'_{k_0}(s)\dot{\gamma}(\bar \chi_{k_0})}{\sigma(s)}=1,\]
where the dot sign denotes the derivative taken respect to $\chi$.  We now proceed to determine the sign of $\sigma(s)$ for $s\in(-\delta,\delta)$.  To this end, we differentiate (\ref{313}) with respect to $\chi$ and then take $\chi=\bar \chi_{k_0}$, then we have that
\begin{equation}\label{314}
\left\{
\begin{array}{ll}
D_1 \Delta \dot u -\phi(\bar u,\bar v) \Delta \Phi_{k_0}-\bar \chi_{k_0} \phi(\bar u,\bar v)\Delta \dot v-\mu\bar u \dot u=\dot \gamma(\bar \chi_{k_0})Q_{k_0}\Phi_{k_0} ,&x\in \Omega,\\
D_2\Delta \dot v-\alpha \dot v+f'(\bar u) \dot u=\dot \gamma(\bar \chi_{k_0}) \Phi_{k_0},&x\in \Omega,
\end{array}
\right.
\end{equation}
where $\dot u=\frac{\partial u(\chi,x)}{\partial\chi}\vert_{\chi=\bar \chi_{k_0}}$ and $\dot v=\frac{\partial v(\chi,x)}{\partial\chi}\vert_{\chi=\bar \chi_{k_0}}$.  Now we multiply (\ref{314}) by $\Phi_{k_0}$ and integrate it over $\Omega$ by parts, then it follows from straightforward calculations and the fact $\int_\Omega \Phi^2_{k_0}dx=1$ that
\begin{equation}\label{315}
\begin{pmatrix}
-D_1 \lambda_{k_0}-\mu\bar{u} & \lambda_{k_0}\bar \chi_{{k_0}}\phi(\bar{u},\bar{v})  \\
~~\\
f'(\bar u) &-D_2\lambda_{k_0}-\alpha
\end{pmatrix}
\begin{pmatrix}
\int_\Omega \dot u \Phi_{k_0}dx \\
~~\\
\int_\Omega \dot v \Phi_{k_0}dx
\end{pmatrix}=\begin{pmatrix}
\dot \gamma(\bar \chi_{k_0})Q_{k_0}-\lambda_{k_0}\phi(\bar u,\bar v)\\
~~\\
\dot\gamma(\bar \chi_{k_0})
\end{pmatrix}.
\end{equation}
We recall from (\ref{210}) that the coefficient matrix in (\ref{315}) in singular.  This implies that
\[\frac{-D_1 \lambda_{k_0}-\mu\bar{u}}{f'(\bar u)}=\frac{\dot \gamma(\bar \chi_{k_0})Q_{k_0}-\lambda_{k_0}\phi(\bar u,\bar v)}{\dot\gamma(\bar \chi_{k_0})}\]\
and consequently
\[\dot\gamma(\bar \chi_{k_0})=\frac{\lambda f'(\bar u)\phi(\bar u,\bar v)}{(D_1+D_2)\lambda +\mu\bar u+\alpha}>0,\]
there it follows from Theorem 1.6 in [Crandall \& Rabinowitz, 1973] and the discussions above that
\begin{equation}\label{316}
\text{sgn} (\sigma(s))=\text{sgn}(-\chi'_{k_0}(0)),
\end{equation}
if $\chi'_{k_0}(0)\neq0$.  When $\chi'_{k_0}(0)=0$, we can easily have from (\ref{316}) and L$'$H$\hat{\text{o}}$pital's rule that $\text{sgn} (\sigma(s))=\text{sgn}(-\chi''_{k_0}(0))$, hence Theorem \ref{theorem31} readily follows from the analysis above.
\end{proof}
The bifurcation branches are transcritical if $ \int_\Omega \Phi_k^3 dx\neq 0$ or $\chi'_k(0)\neq 0$.  If domain $\Omega$ is an interval or multi--dimensional rectangle, the eigenfunction $\Phi_k$ is cosine or a product of cosine functions, therefore $\int_\Omega \Phi_k^3 dx=0$ and $\chi'_k(0)=0$, and the bifurcation branch is of pitchfork type, i.e., being one sided.

\section{Numerical simulations in 2D}\label{section4}
This section is devoted to the numerical studies of development of spatially inhomogeneous solutions to the time--dependent system (\ref{11}) over 2D domains.  In particular, we choose $\Omega=(0,L)\times (0,L)$.  In this special case, the Neumann eigenvalue problem (\ref{21}) has eigen--pairs
\[\lambda_{mn}=(m^2+n^2)\pi^2/L^2,\Phi_{mn}(x,y)=\cos \frac{m\pi x}{L} \cos \frac{n\pi y}{L},~m,n\in\mathbb N^+.\]
Each eigenvalue gives rise to bifurcation point $\chi_{mn}$, with $\lambda_k$ being replaced by $\lambda_{mn}$ in (\ref{210}).  Though the boundary of the square is not smooth at the corners, one can interpret the boundary condition in the weak sense via Green's identity--see the discussions in [Chertock, \emph{et al}., 2012].  We also want to mention that numerical studied in 1D model with volume filling effect has been performed in [Ma, \emph{et al}., 2012].
\begin{figure}[!htb]\centering
  \includegraphics[width=3.8in]{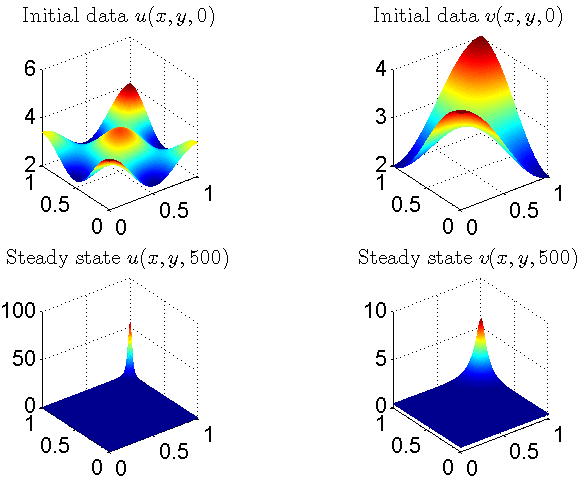}
 \caption{Emergence of single boundary spike at corner $(1,1)$.}\label{fig2}
\end{figure}
We are interested in patterns with concentrating properties in contrast to the small amplitude bifurcating solutions established in Theorem \ref{theorem21}.  Our numerical studies suggest that for properly chosen parameters and initial data, the stable steady states can develop quite complicated structures such as spikes, stripes, etc.  These patterns with concentrating structures can be realistic modelings of cellular aggregations, though rigorous analysis is needed to fully understand the initiation, development and evolution of these striking structures.  In particular, in order to illustrate the effect of system parameters $D_1$, $D_2$, $\chi$ and the initial data on the pattern formations, we restrict $\alpha=1$ and $f(u)=u$ in this section.

First of all, we numerically study (\ref{11}) that develops a boundary spike at the corner $(1,1)$ of $\Omega=(0,1)\times (0,1)$ in Figure \ref{fig2}.  We choose $D_1=\chi=5$, $D_2=0.01$, $\bar u=3$ and $\mu=1$; the initial data are chosen to be $u_0=\bar u+\cos 2\pi x\cos \pi y$ and $v_0=\bar v+\cos \pi x\cos \pi y$.  Our numerical studies suggest that small chemical diffusion rate $D_2$ supports steady states with spikes.  Numerical results in Figure \ref{fig2} also suggest that spikes are usually expected to develop and stabilize at the location where initial data maximize.  This is also supported by results in Figures \ref{fig4} and Figure \ref{fig5}.

\begin{figure}[!htb]\centering
  \includegraphics[width=3.8in]{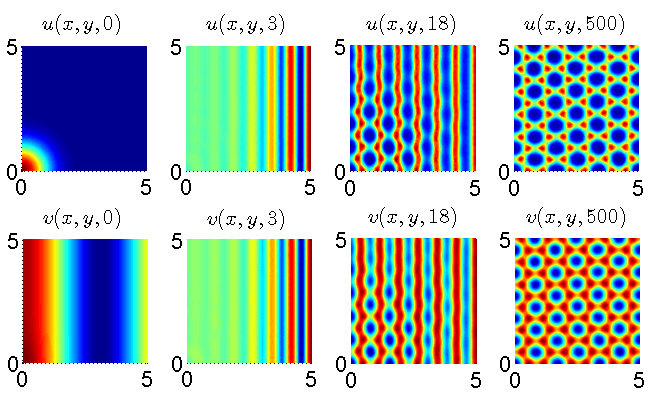}
 \caption{Emergence of hexagon steady states through unstable stripes.}\label{fig3}
\end{figure}

\begin{figure}[!htb]\centering
  \includegraphics[width=3.8in]{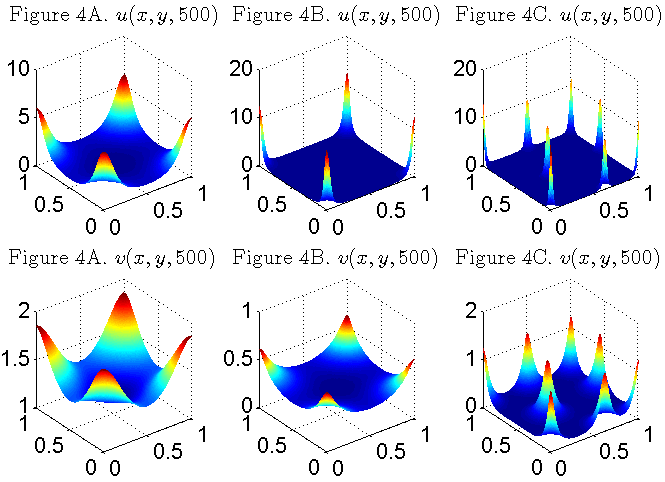}
 \caption{Stable steady states with boundary/interior spikes with the same initial data, but different sets of system parameters.}\label{fig4}
\end{figure}

\begin{figure}[!htb]\centering
  \includegraphics[width=4.2in]{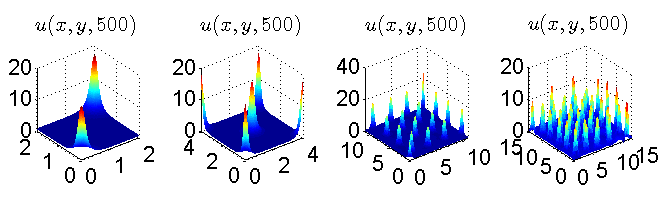}
 \caption{Stable steady states with boundary and interior spikes as domain size increases, with parameters and initial data being fixed.  This experiment indicates that large domain supports the emergence of spikes.}\label{fig5}
\end{figure}

\begin{figure}[!htb]\centering
  \includegraphics[width=4.2in]{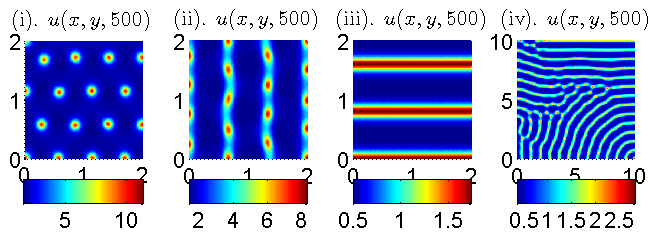}
 \caption{Stable steady states with spikes, dotted--stripes, stripes etc. emerge.}\label{fig6}
\end{figure}

In Figure \ref{fig3}, we study (\ref{11}) with cell diffusivity being $D_2$ relatively small with $\Omega=(0,1)\times(0,1)$.  In particular, we choose $D_1=0.0625$, $D_2=1$, $\chi=19$, $\mu=8$ and $\bar u=1$; the initial data are taken to be $u_0=\bar u+e^{-(x^2+y^2)}$ and $v_0=2\bar u+0.01\cos x$.  It is presented in Figure \ref{fig3} that stripes develop at $t\approx3$ and break down at $t\approx18$; then the stripes lose stability to steady state with hexagon blocks, which stabilizes for $t$ being large enough.  We want to mention that the hexagon steady states are also numerically presented in [Osaki, \emph{et al}., 2002], which motivates our results in Figure \ref{fig3}.
In Figure \ref{fig4}, we investigate the effect of system parameters on the formations of stable steady states.  For this purpose, we fix $D_1=1$ and $\mu=10$ and choose the initial data to be small spatial perturbations from $(\bar u,\bar v)=(3,3)$, then we perform extensive numerical studies of (\ref{11}) as $D_2$ and $\chi$ vary.  In Figures 4A, 4B and 4C, $(D_2, \chi)$ are selected to be $(0.1,5)$, $(0.1, 20)$ and $(0.005,5)$ respectively.  We observe that both small $D_2$ and large $\chi$ create spikes; Figure 4A and Figure 4B suggest that large $\chi$ drives spikes to concentrate at corners;  Figure 4A and Figure 4C indicate that small $D_1$ supports the emergence of stable boundary and interior spikes.  These results also indicate that large chemotactic rate $\chi$ shrinks the region where the spike is supported.  It is quite reasonable to predict that the spike concentrates at a single point or multiple points as $\chi \rightarrow \infty$, however our numerical experiments suggest that blow--up may occur in (\ref{11}) if $\chi$ is sufficiently large.  Rigorous proof is needed to support this observation.  All the numerical results demonstrate that small spatial perturbations of homogeneous equilibrium can lead to stable patterns with various complicated and interesting structures.
We observe that, by reflecting and periodically extending steady state with a single spike, one can construct multi--spike solutions thanks to the homogeneous Neumann boundary conditions.  This motivates us to investigate the formations of stable steady states in larger square domains.  In Figure \ref{fig5}, we plot stable spiky steady states as the domain size $L$ increases.  The parameters are chosen to be $D_1=\chi=5$, $D_2=0.1$, $\mu=1$ and $\bar u=\bar v=3$; initial data are $u_0=\bar u+\cos (2x+1)\cos (2y+1)$ and $v_0=\bar v+\cos (2x) \cos (2y)$.  For $\Omega=(0,2)\times (0,2)$, two spikes emerge on the opposite diagonal corners; for $\Omega=(0,4)\times (0,4)$, three more stable spikes develop, with two on the corners and one in the center;  for $\Omega=(0,10) \times(0,10)$ and $(0,15)\times (0,15)$, we see more interior spikes emerge over the habitat.  It is realistic to claim that infinitely many spikes can emerge if $\Omega$ expands to the whole domain, however, rigorous analysis of the mechanism that drives the formations of spikes is an interesting but also challenging question that deserves future explorations.
Finally, we present some stable numerical solutions with multi--spikes and/or multi--stripes in Figure \ref{fig6}.  These patterns are heuristic in understanding the dynamics of system (\ref{11}) with logistic growth, demonstrating that it can develop many other complicated and interesting structures besides those established in previous numerical experiments.  The numerical results and choice of system parameters in this figure are motivated by those obtained in [Kuto, \emph{et al}., 2012].  For all plots in Figure \ref{fig6}, the initial data are taken to be $u_0=\bar u+0.05 e^{-\frac{x^2+y^2}{2}}$, $v_0=\bar v+0.05 e^{-\frac{x^2+y^2}{2}}$.  In (i), we choose that $D_1=D_2=0.25$, $\mu=5$, $\bar u=3$, $\chi=10$; stable steady solution emerge with multi--spikes, boundary and interior.  In (ii), we slightly change the diffusion rates to be $D_1=0.125$ and $D_2=0.5$ and keep the rest parameters the same as those in (i), then stable multi--spikes also develop, which form a stable strip--like patterns.  This motivates us to set $D_1$ smaller in plot (iii).  We set $D_1=0.0625$, $D_2=1$, $\mu=6$, $\bar u=1$ in plot (iii), then we see that regularized stripes form a stable pattern.  Plot (iv) is obtained numerically by enlarging the domain size in plot (iii), while the parameters in both plots are set to be the same.

\section{Conclusion and discussion}\label{section5}
In this paper, we study positive steady states of Keller--Segel chemotaxis system (\ref{12}) with logistic growth over multi--dimensional bounded domain.  We show that the homogeneous solution $(\bar u,\bar v)$ loses its stability as chemo--attraction rate $\chi$ surpasses a threshold value $\chi_0$.  Local bifurcation analysis is performed to establish the existence of nonconstant positive solutions $(u_k(s,x),v_k(s,x))$ of (\ref{12}) in Theorem \ref{theorem21}.   It is shown that the first steady state bifurcation occurs at the same location where $(\bar u,\bar v)$ loses its stability.  Our results indicate that chemotaxis is a leading mechanism that drives the formation of these small amplitude bifurcating solutions.

We also investigate stability of the bifurcating solutions around $(\bar u,\bar v,\chi_k)$ in Theorem \ref{theorem31}.  Our result states that the only stable bifurcation branch must have a wavemode number $k_0$ which minimizes the bifurcation value $\bar \chi_k$ in (\ref{210}) over $\mathbb N^+$.  This provides a selection mechanism of stable wavemode to which the homogeneous solution $(\bar u,\bar v)$ loses its stability.  The same results have been obtained in [Ma, \emph{et al}. , 2012] for chemotaxis system with volume--filling effect over 1D (though it is easy to see that their results carry over for more general systems).  On one hand, our existence and stability results extend theirs to multi--dimension.  On the other hand, we show that the sign of $\int_\Omega \Phi_k^3dx$ determines the shape of the transcritical bifurcation branches if the integral is not zero.  If the domain has a geometry such that $\int_\Omega \Phi_k^3dx=0$, for example when it is an interval or is a multi--dimensional rectangle, then the bifurcation branch is of pitch--fork type in this case.  Furthermore, we evaluate $\chi_k''(0)$ in terms of several integrals which can be computed involving system parameters.  Finally, we present numerical simulations of formation of patterns with stable structure such as with boundary spikes, interior spikes, stripes, etc.  Our computational results also suggest that the logistic growth chemotaxis system (\ref{11}) can develop many striking structures which can model cellular aggregations that emerge in realistic biological systems.

There are also some interesting questions that deserve future explorations on system (\ref{11}) and its stationary counterpart (\ref{12}).  The appearance of the logistic growth is an important mechanism to inhibit solutions from blowing up in finite or infinite time, however, whether or not this is sufficient remains an open problem, in particular when both $\chi$ and the space dimension $N$ are large.  Further questions can be asked about the global existence of (\ref{11}) with cellular kinetics different from the logistic type.  According to the global bifurcation theory in [Crandall \& Rabinowitz, 1971] and [Shi \& Wang, 2009], the continuum of each bifurcation branch satisfies one of three alternatives, and rigorous characterization is very important for the analysis of positive solutions of elliptic systems, however, highly nontrivial mathematical techniques are required to deal with this problem for multi--dimensional domains.

Our numerics illustrate the emergence and development of spatial patterns with spiky and other concentrating properties.  It is suggested in literature that both large chemo--attraction rate $\chi$ and small chemical diffusion rate $D_2$ support the emergence and stability of these nontrivial patterns.  Rigorous construction and/or mathematical analysis of solutions with concentrating structures is an interesting and challenging problem to work on, even over one--dimensional finite domain.  Moreover, rigorous stability analysis of these patterns is apparently another delicate problem that deserves probing in the future.

\end{document}